\documentclass[aop,MSNbibl,seceqn,dvips]{arximspdf}
\usepackage{graphicx}
\doi{10.1214/13-AOP841} 
\volume{42}
\issue{3}
\pubyear{2014}
\firstpage{1285}
\lastpage{1296}

\makeatletter
\newcommand{\eqref}[1]{(\ref{#1})}
\newcommand{\N}{{\mathbb N}}
\newcommand{\Z}{{\mathbb Z}}
\newtheorem{theorem}{Theorem}[section]
\newtheorem{proposition}{Proposition}[section]
\newtheorem{lemma}{Lemma}[section]
\newproclaim{remark}{Remark}[section]
\makeatother

\begin{document}
\begin{frontmatter}

\title{Convergence in distribution for subcritical 2D~oriented
percolation seen from its rightmost~point}
\runtitle{Oriented percolation}

\begin{aug}
\author{\fnms{E.~D.} \snm{Andjel}\corref{}\thanksref{t1}\ead[label=e1]{andjel@impa.br}}
\thankstext{t1}{Visiting IMPA. Supported by CAPES.}
\runauthor{E.~D. Andjel}
\affiliation{Universit\'e d'Aix-Marseille}
\address{Universit\'e d'Aix-Marseille\\
LATP, 39 Rue Joliot-Curie\\
13453 Marseille, cedex 13\\
France\\
\printead{e1}}
\end{aug}

\received{\smonth{5} \syear{2012}}
\revised{\smonth{2} \syear{2013}}

%
\begin{abstract}
We study subcritical two-dimensional oriented percolation seen from its
rightmost point on the set of infinite configurations
which are bounded above.
This a Feller process whose state space is not compact and has no
invariant measures.
We prove that it converges in distribution to a measure which charges
only finite configurations.
\end{abstract}

%


\end{frontmatter}

\section{Introduction and main results}\label{sec:intro}
\subsection{Background}\label{sec1.1}
Two-dimensional oriented percolation and its continuous time analog the
one-dimensional contact process, seen from their rightmost point, have
been studied in several papers.
Durrett \cite{D} proved that in the critical and supercritical phase
there exists an invariant measure. Then, Schonmann proved that there
are no such measures
in the subcritical phase \cite{S}.
These two papers consider only the discrete time model, but their
results hold also for some continuous time one-dimensional process
which include the contact process
(see \cite{ASS}). Galves and Presutti \cite{GP} proved that the
one-dimensional contact process seen from the rightmost point converges
in the supercritical phase to a unique invariant measure.
This last result was then
extended by Cox, Durrett and Schinazi \cite{CDS}
to the critical phase. There are no difficulties in adapting these
convergence results to the discrete time setting. Finally, we mention
\cite{GP}
and \cite{K} where the position of the rightmost point is shown to
satisfy a central limit theorem. In this paper, we prove that
convergence of the discrete time process seen from the rightmost point
also occurs in the subcritical phase although there are no invariant measures.

\subsection{Definitions}\label{sec1.2}

Let
%
%
\begin{equation}
\Lambda=\bigl\{(x,y)\dvtx x, y \in\Z, y\geq0, x+y \in2\Z\bigr\}.\vadjust{\goodbreak}
\end{equation}

Draw oriented bonds from each point $(m,n)$ in $\Lambda$ to
$(m+1,n+1)$ and
to $(m-1,n+1)$. In this paper, we suppose that bonds
are open independently of each other and that each bond is open with
probability $p\in(0,1)$. To formalise this, let~$\mathcal{B} $ be the
set of all bonds with both endpoints in $\Lambda$ and assume that
$\{Z_b\dvtx b \in\mathcal{B}\} $ is a collection of i.i.d. random
variables whose distribution is a Bernoulli with parameter $p$. A bond
$b$ will be considered open (closed) if
$Z_b=1$ ($Z_b=0$). The event consisting on the existence of an open
path from $A$ to $B$, where $A$ and $B$ are subsets of $\Lambda$, will
be denoted by
$\{A\rightarrow B\}$ and its complement by $\{A\nrightarrow B\}$. When
either $A$ or $B$ or both are singletons, say $\{x\}$ and $\{y\}$,
respectively, we will write
$\{x\rightarrow y\}$, $\{x\rightarrow B\}$, etc.

Given a subset $A$ of $2\Z$ we let
%
%
\begin{eqnarray}
\xi^A_n= \bigl\{y\dvtx(y,n)\in\Lambda\mbox{ and } (x,0)
\rightarrow(y,n) \mbox{ for some } x\in A\bigr\},
\nonumber
\\[-8pt]
\\[-8pt]
\eqntext{n=0,1,\ldots.}
\end{eqnarray}
Then, $(\xi^A_n, n\geq0) $ is a Markov chain taking values in the
subsets of $2\Z$ at even times and of $2\Z+1$ at odd times.

Let $A$ be an infinite subset of $2\Z$ such that $\sup A<\infty$.
Then, for all $n>0$, the supremum of $\xi^A_n$ is finite and a simple
Borel--Cantelli argument
shows that $\xi^A_n$ is a.s. infinite. For such initial conditions, we let
\[
r\bigl(\xi^A_n\bigr) =\sup\xi^A_n
\quad\mbox{and}\quad \zeta^A_n= \bigl\{x-r\bigl(
\xi^A_n\bigr)\dvtx x \in\xi^A_n
\bigr\}.
\]
Then $(\zeta^A_n, n\geq0)$ is a Markov chain on infinite subsets of
$2\Z_-=:\{0,-2,\break -4,\ldots\}$
containing $0$. For finite subsets $A$ we may also define the Markov
chain $(\zeta^A_n, n\geq0)$ by simply adopting the convention: $\zeta
^A_n=\varnothing$ if $\xi^A_n=\varnothing$. Obviously, $\varnothing
$ is an
absorbing state for both
$(\xi^A_n, n\geq0) $ and $(\zeta^A_n, n\geq0) $.

In the sequel,
%
%
\begin{eqnarray}
S&=&\{ \mbox{infinite subsets of } 2\Z_- \mbox{ containing } 0\},
\\
S_0&=&\{ \mbox{finite subsets of } 2\Z_- \mbox{ containing } 0\}
\end{eqnarray}
and
%
%
\begin{equation}
\bar S=\{ \mbox{subsets of } 2\Z_- \}.
\end{equation}

We will consider $S$ and $S_0$ as subsets of $\bar S$ which we identify
with $\{0,1\}^{2\Z_-}$ by means of the bijection:
$F(A)={\bf1}_A$. Then, $\bar S$ inherits the product topology of $\{
0,1\}^{2\Z_-}$ and becomes a compact space. The subsets $S$ and $S_0$
of $\bar S$ are now endowed with the induced topology. Probability
measures on either $S$ or $S_0$ will be seen as measures on
$\bar S$ and the space of all probability measures on $\bar S$ will be
endowed with the topology of weak convergence.

Standard coupling arguments show that $P(\xi^0_n\neq\varnothing\mbox{ for all }n) $ increases with $p$ and we can define the critical value
$p_c$ of the parameter $p$ as
the supremum of its values for which the above probability is $0$. It
is well known (see \cite{D}) that $0<p_c<1$.
Throughout this paper, we assume that $p\in(0,p_c)$.

\subsection{Theorems}\label{sec1.3}
Before stating our results, we recall that a quasi-stationary
distribution of a Markov chain $(X_n;n\geq0)$ on $S_0\cup\{
\varnothing\}
$ with absorbing
state $\varnothing$ is a probability measure $\nu$ on $S_0$ such that
$P_{\nu}(X_n=x\vert T>n)=\nu(x)$ for all $n\in\N$ and $x\in S_0$,
where $T=\inf\{k\dvtx X_k=\varnothing\}$. We
refer the reader to \cite{FKM} for more information concerning
quasi-stationary distributions.
Our first theorem is not new, it is immediately obtained from Theorem 1
of \cite{FKM}.

%
\begin{theorem}\label{T2} Suppose $0<p<p_c$ and let $T=\inf\{n\dvtx
\zeta
^0_n=\varnothing\}$. Then the conditional distribution
of $\zeta^0_n$ given $\{T>n\}$ converges as n goes to infinity to a
probability measure $\nu$ on $S_0$. Moreover, $\nu$ is the minimal
quasi-stationary distribution of the $\zeta_n$ process on
$S_0\cup\{\varnothing\}$.
\end{theorem}

We now state our main result which was conjectured by Galves, Keane and
Meilijson.
As expected by the authors of \cite{FKM} (see Remark 7 in page 606 of
that reference), Theorem \ref{T2} is the key ingredient to prove it.

%
\begin{theorem}\label{T1} Suppose $0<p<p_c$. Then, for any $A\in S$ the
distribution of $\zeta^A_n$ converges as n goes to infinity to $\nu$,
where $\nu$ is as in Theorem \ref{T2}.
\end{theorem}

The paper is organised as follows: Section \ref{sec2} starts explaining the
strategy we will follow, continues stating two lemmas and then deduces
Theorem \ref{T1} from these lemmas.
Then, in Section \ref{sec3} we prove those two lemmas.
\section{\texorpdfstring{Proof of Theorem \protect\ref{T1}}
{Proof of Theorem 1.2}}\label{sec2}

We start this section introducing some more notation:
Let $f$ be real-valued function defined on $S\cup S_0$.
We say that $f$ is a cylinder function depending only on coordinates
$-2r,\ldots,-2$ if there exists a function $g$ defined on subsets of
$\{-2r,\ldots,-2\}$ such that $f(A)=g(A\cap\{-2r,\ldots,-2\}) $ for all
$A\in S\cup S_0$. For such functions, $\| \|$ will denote
the supremum norm
\[
\|f\|=\sup_{A\in S\cup S_0}\bigl\vert f(A)\bigr\vert.
\]
For $(x,m)\in\Lambda$ let
\[
C_{x,m}= \bigl\{(y,k)\in\Lambda\dvtx k\geq m, \vert y-x\vert\leq k-m
\bigr\}
\]
and call this set the cone emerging from $(x,m)$.
For $r\in\N$,
call level $r$ the set
\[
L_r=\bigl\{(x,n)\in\Lambda\dvtx n=r \bigr\}.
\]
We will say that level $n$ is higher than level $m$ if $n\geq m$.
In the sequel, $\nu_r$ will be the distribution of $\zeta^0_r$ given
$\{
T>r\}$ where $T=\inf\{k\dvtx\zeta_k^0=\varnothing\}$ and
$A$ will be fixed but arbitrary element in $S$.

We now sketch the proof of Theorem \ref{T1}:
We first find the rightmost point $x_0$ of $A$ satisfying $\xi
^{x_0}_n\neq\varnothing$. We would like to apply Theorem \ref{T2} but
cannot do it immediately because
we are conditioning not only on $\{\xi^{x_0}_n\neq\varnothing\}$ but
also on $\{\xi^{y}_n= \varnothing\}$ for all $y\in A\cap\{z\dvtx
z>x_0\}$.
However since $p<p_c$, there is a positive
probability that no point of $C_{x_0,0}$ can be attained from $\{
(y,0)\dvtx y>x, y\in A\}$ following open paths. If this occurs, then the
distribution of $\zeta^{x_0}_n$ is $\nu_n$. If this
fails to happen
we look at the highest level in $C_{x_0,0}$ attained from $\{(y,0)\dvtx
y>x_0, y\in A\}$ and repeat the argument from that level. We keep doing
so until the corresponding emerging cone is not attained.
Once this happens, we will derive from $p<p_c$ that the elements of
$\xi
_n^A\setminus\xi^{x_0}_n$ are far to the left of $\xi^{x_0}_n$ for
large $n$.
In carrying out this approach, the main difficulty
comes from keeping track of several conditionings.
To make this argument rigorous, we begin defining two
sequences of r.v.'s $Y_i$ and $X_i$
as follows: Let
%
%
\begin{equation}
\label{a} Y_0=0
\end{equation}
and
%
%
\begin{equation}
\label{aa} X_0=\sup\bigl\{x\in A\dvtx(x,0)\rightarrow
L_n\bigr\}.
\end{equation}
Then, given $Y_0,Y_1,\ldots,Y_i$ and $X_0,\ldots,X_i$, we let
%
%
\begin{eqnarray}
\label{aaa}&& Y_{i+1} = \sup\bigl\{k\dvtx\exists u,v\dvtx
(u,Y_i)\rightarrow(v,k)
\nonumber
\\[-8pt]
\\[-8pt]
\nonumber
&&\hspace*{54pt}{}\mbox{with } u>X_i,
u\in\xi_{Y_i}^A \mbox{ and } (v,k)\in C_{X_i,Y_i}
\bigr\}\vee Y_i
\end{eqnarray}
and
%
%
\begin{equation}
\label{aaaa} X_{i+1}=\sup\bigl\{x\dvtx x\in\xi_{Y_{i+1}}^A
\mbox{ and } (x,Y_{i+1})\rightarrow L_n\bigr\}.
\end{equation}
%
Note that $Y_i\leq Y_{i+1}$ and that as soon as $Y_i=Y_{i-1}$ both
sequences become constant. The reader may find helpful to
have now a first look at Figure \ref{fig1} in the next section.

We now state two lemmas which will be proved in the next section. In
the second of these lemmas, we use the fact that
on the event $\{0\rightarrow L_n\}$ the process $(\zeta^0_k,
k=0,\ldots
,n)$ takes values on $S_0$.
%
%
\begin{lemma}\label{c1}
Let $I=\inf\{i\dvtx Y_i=Y_{i+1}\}$. Then, there exists $\beta>0$ such that
for all $m$:
\begin{longlist}[(a)]
\item[(a)] $P(I\geq m)\leq\exp(-\beta m)$,

\item[(b)] $P(Y_I\geq m^2) \leq(m+1) \exp(-\beta m) $ and

\item[(c)] $P(I\leq m, Y_I\leq m^2) \geq1-(m+2) \exp(-\beta m) $.
\end{longlist}
\end{lemma}

%
\begin{lemma}\label{c6}
Let $A$ be an element of $S$, let $f$ be a cylinder function on $S\cup
S_0$ depending only on the coordinates $-2r,\ldots,-2$ and let $\beta$
be as in Lemma \ref{c1}. Then, for all $i,j\leq n$,
\[
\bigl\vert E\bigl(f\bigl(\zeta_n^A\bigr)\vert I=i,
Y_i=j\bigr)-E\bigl(f\bigl(\zeta^0_n\bigr)
\vert0\rightarrow L_n\bigr)\bigr\vert\leq2\|f\| (n+r)\exp\bigl(-\beta(n-j)
\bigr).
\]
\end{lemma}

We now proceed to prove our main result.\vadjust{\goodbreak}

\begin{pf*}{Proof of Theorem \ref{T1}}
Let $f$ be a cylinder function on $S\cup S_0$ depending only on the
coordinates $-2r,\ldots,-2$ and let $m=\lfloor n^{1/3}\rfloor$ where
$\lfloor\cdot\rfloor$ denotes the integer part.
By part (c) of Lemma \ref{c1}, we have
\begin{eqnarray*}
&&\Biggl\vert E\bigl(f\bigl(\zeta_n^A\bigr)\bigr)-\sum
_{i=0}^m \sum_{j=0}^{m^2}E
\bigl(f\bigl(\zeta_n^A\bigr)\vert I=i,Y_i=j
\bigr)P( I=i,Y_i=j)\Biggr\vert
\\
&&\qquad\leq\|f\|(m+2)\exp(-\beta m).
\end{eqnarray*}
%

Therefore,
\begin{eqnarray*}
&&\bigl\vert Ef\bigl(\zeta_n^A\bigr)-E\bigl(f\bigl(
\zeta^0_n\bigr)\vert0\rightarrow L_n\bigr)
\bigr\vert
\\
&&\qquad\leq\sum_{i=0}^m \sum
_{j=0}^{m^2} \bigl( \bigl\vert E\bigl(f\bigl(
\zeta_n^A\bigr)\vert I=i,Y_i=j\bigr)-E\bigl(f
\bigl(\zeta^0_n\bigr) \vert0\rightarrow L_n
\bigr) \bigr\vert \bigr) P(I=i,Y_i=j)\\
&&\qquad\quad{} +
2\|f\|\bigl(1-P\bigl(I\leq m,Y_I\leq m^2\bigr)\bigr)
\\
&&\qquad\leq\sum_{i=0}^m \sum
_{j=0}^{m^2} \bigl( 2\| f\|(n+r)\exp\bigl(-\beta(n-j)
\bigr)P(I=i,Y_i=j) \bigr)\\
&&\qquad\quad{}+
2\|f\|\bigl(1-P\bigl(I\leq m,Y_I\leq m^2\bigr)\bigr)
\\
&&\qquad\leq2\|f\|(n+r)\exp\bigl(-\beta\bigl(n-m^2\bigr)\bigr)+2\|f\|\bigl(1-P\bigl(I
\leq m,Y_I\leq m^2\bigr)\bigr),
\end{eqnarray*}
where the second inequality follows from Lemma \ref{c6}.
Since $m=\lfloor n^{1/3}\rfloor$ this and part (c) of Lemma \ref{c1},
imply that
\[
\lim_n \bigl\vert Ef\bigl(\zeta_n^A
\bigr)-E\bigl(f\bigl(\zeta^0_n\bigr)\vert0\rightarrow
L_n\bigr)\bigr\vert=0,
\]
and the result follows from Theorem \ref{T2}.
\end{pf*}

\section{\texorpdfstring{Proofs of Lemmas \protect\ref{c1} and \protect\ref{c6}}
{Proofs of Lemmas 2.1 and 2.2}}\label{sec3}

In this section, for $r\in\N$ and $x\in2\Z_-$, ${\mathcal G}_x^r$
will denote the $\sigma$-algebra generated by the random variables
which determine the state of the bonds with both vertices
in $(\bigcup_{i=0}^{-x/2}C_{x+2i,0})\cap(\bigcup_{j=0}^r L_j)$,
${\mathcal
G}^r$ will denote the $\sigma$-algebra generated by the random
variables which determine the state of the bonds with both vertices
in $\bigcup_{i=0}^r L_i$ and ${\mathcal G'}^r$ will denote the $\sigma
$-algebra generated by the random variables which determine the state
of the bonds with both vertices
in $\bigcup_{i=r}^{\infty} L_i$.
Besides this, an event belonging to a $\sigma$-algebra generated by
random variables determining the state of a finite number of bonds will
be called an \textit{elementary cylinder} of that
$ \sigma$-algebra if it is nonempty and does not contain any nonempty
proper subset of that $\sigma$-algebra.
The first lemma of this section is an immediate consequence of the
exponential decay of $P(\xi^x_n\neq\varnothing)$ (see Section 7 of
\cite
{D}) and we omit its proof.
%
%
\begin{lemma}\label{l0}
There exists a constant $\beta>0$ such that for all $x\in2\Z$ and all
$m\in\N$ we have
\[
P \Biggl((y,0)\rightarrow C_{x,0}\cap\Biggl(\bigcup
_{j=m}^{\infty} L_j\Biggr)\mbox{ for some
}y>x \Biggr)\leq\exp(-\beta m).
\]
\end{lemma}
%

In our next lemma, for notational convenience we let $r_0=0$
and recalling \eqref{a}--\eqref{aaaa}, consider events of the form
\begin{eqnarray*}
G(x_0)&=&\{X_0=x_0\} \qquad\mbox{and for } i\geq1
\\
G(r_1,\ldots,r_i;x_0,\ldots,x_i)&=&
\{Y_1=r_1,\ldots, Y_i=r_i,
X_0=x_0,\ldots,X_i=x_i\},
\end{eqnarray*}
where $0\leq r_1\leq\cdots\leq r_i$ are integers and
$(x_0,0),(x_1,r_1),\ldots, (x_i,r_i)\in\Lambda$.

%
\begin{lemma}\label{l10} Let $i$ be a nonnegative integer and
let $F$ be an elementary cylinder in ${\mathcal G}_{x_0}^{r_i}$ having
a nonempty intersection with $G(r_1,\ldots,r_i;x_0,\ldots,x_i)$.
Then, $L_{r_i}$ contains $i+2$
finite subsets $A_{i,1},\ldots, A_{i,i},B_i,D_i$ determined by
$F,n,x_0,\ldots,x_i$, $r_1,\ldots,r_i$ only and such that
%
%
\begin{eqnarray}
\label{eq20} %
&& F\cap G(r_1,\ldots,r_i;x_0,
\ldots,x_i)\nonumber\\
&&\qquad=F\cap\bigl\{(x_i,r_i)\rightarrow
L_n\bigr\} \cap\{u\nrightarrow L_n\ \forall u\in
D_i\}\\
&&\qquad\quad{} \cap\{u \nrightarrow L_{r_i+1}\ \forall u\in
B_i\}\cap\Biggl(\bigcap_{j=1}^i
\{u\nrightarrow L_n,u\nrightarrow C_{x_{j-1},r_{j-1}}\ \forall u\in
A_{i,j}\} \Biggr).\nonumber
\end{eqnarray}
Moreover,
%
%
\begin{equation}
\label{eq201} \xi_{r_i}^A(x_i)=1, \bigl\{z\dvtx
z>x_i,\xi_{r_i}^A(z)=1\bigr\}=D_i
\cup B_i \cup\Biggl(\bigcup_{j=1}^i
A_{i,j}\Biggr)
\end{equation}
and
%
%
\begin{eqnarray}
\label{eq202} %
x_i<d_i<b_i<a_{i,i}<a_{i,i-1}<
\cdots<a_{i,1}
\nonumber
\\[-8pt]
\\[-8pt]
 \eqntext{\forall d_i\in D_i,b_i
\in B_i, a_{i,j}\in A_{i,j},j=1,\ldots, i.}
\end{eqnarray}
%
\end{lemma}
\begin{remark*} If $F$ is disjoint of $G(r_1,\ldots,r_i;x_0,\ldots,x_i)$,
we may extend the definition of the sets $A_{i,1},\ldots,
A_{i,i},B_i,D_i$ by letting them be the empty set.
In this way, they become random $\mathcal{G}_{x_0}^{r_{i}}$-measurable
sets, hence independent of the
$\sigma$-algebra~${\mathcal G'}^{r_i}$.
\end{remark*}

\begin{pf*}{Proof of Lemma \ref{l10}}
To follow this proof, we recommend the reader to look at Figure \ref{fig1}. This
may help visualizing the different sets involved in it.
We proceed by induction on $i$. If $i=0$, then ${\mathcal
G}_{x_0}^{r_i}={\mathcal G}_{x_0}^{0}$ is the trivial
$\sigma$-algebra. Hence, $F$ must be the whole probability space and
the statement holds with $D_0=\{(z,0)\dvtx z>x_0,z\in A\}$ and
$B_0=\varnothing$.
Assume the statement holds for some given $i$, and let $F'$ be an
elementary cylinder of $\mathcal{G}_{x_0}^{r_{i+1}}$. Call $F$ the
unique elementary cylinder of
$\mathcal{G}_{x_0}^{r_{i}}$ which contains $F'$. Then, by the inductive
hypothesis there are $i+2$ subsets $A_{i,1},\ldots, A_{i,i},B_i, D_i$ for
which (\ref{eq20}),
(\ref{eq201}) and (\ref{eq202}) hold.
Now, define the following subsets of $L_{r_{i+1}}$:
\begin{eqnarray*}
A_{i+1,j}&=&\bigl\{(x,r_{i+1})\dvtx A_{i,j}\rightarrow(x,r_{i+1}) \bigr\}\qquad
(i\geq1,j=1,\ldots,i),
\\
A_{i+1,i+1}&=&\bigl\{(x,r_{i+1})\notin C_{x_i,r_i}\dvtx D_i \rightarrow
(x,r_{i+1})\bigr\} \qquad(i\geq0),\\
B_{i+1}&=& \bigl\{(x,r_{i+1})\in C_{x_i,r_i}\dvtx D_i \rightarrow
(x,r_{i+1})\bigr\}\quad
\mbox{and}\\
D_{i+1}&=& \bigl\{(x,r_{i+1})\in C_{x_i,r_i}\dvtx x>x_{i+1},(x_i,r_i)
\rightarrow
(x,r_{i+1})\bigr\}\setminus B_{i+1}.
\end{eqnarray*}

%
\begin{figure}

\includegraphics{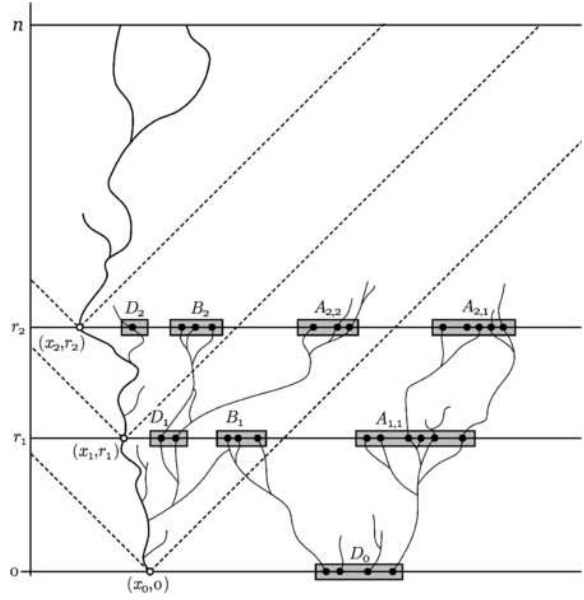}

\caption{Here
$I=2$, $(X_0,Y_0)=(x_0,0)$, $(X_1,Y_1)=(x_1,r_1)$, $(X_2,Y_2)=(x_2,r_2)$,
the doted lines are the emerging cones from these three points and the
full lines are the open paths starting from $A$.}
\label{fig1}
\end{figure}

It is now tedious but straightforward to verify that these
sets satisfy (\ref{eq20}),
(\ref{eq201}) and (\ref{eq202}) with $F'$ and $i+1$ replacing $F$ and
$i$, respectively.
\end{pf*}
%

%
%
\begin{proposition}\label{p1}
Let $\beta$ be as in Lemma \ref{l0}. Then, for all $m,i$ and all
$x_0,\ldots, x_i,r_1,\ldots,r_i$ we have:
%
%
\begin{equation}
\label{eq21} P\bigl(Y_{i+1}-Y_i\geq m\vert
G(r_1,\ldots,r_i,x_0,\ldots,
x_i)\bigr)\leq\exp(-\beta m).
\end{equation}
\end{proposition}

\begin{pf}
Call $\Upsilon$ the set of all paths from $(x_i,r_i)$ to $L_n$.
Given a path $\gamma\in\Upsilon$, call $ A_{\gamma, r}$ the oriented
graph composed by the bonds having both vertices in $\bigcup_{j=r_i}^n
L_j$ and at least one vertex
strictly to the right of $\gamma$, and by the vertices of such bonds.
Let $\Gamma_{r}$ be the rightmost open path
starting from $(x_i,r_i)$ and attaining $L_n$. Let $\gamma$ be a
possible value of $\Gamma_{r}$.
Note that the event $\{\Gamma_{r}=\gamma\}$ is constituted by the
configurations for which $\gamma$ is open but there is no open path
from a vertex of $\gamma$ contained in
$ A_{\gamma, r}$ and reaching either $L_n$ or another point in $\gamma
$. Hence, the event $\{\Gamma_{r}=\gamma\}$ is the intersection of
the event
$\{\gamma\mbox{ is open}\}$ and a decreasing event $D(\gamma)$ on the
graph~$ A_{\gamma, r}$.

Let $F$ be an elementary cylinder in ${\mathcal G}_{x_0}^{r_i}$ having
a nonempty intersection with $G(r_1,\ldots,r_i;x_0,\ldots,x_i)$. To
prove the proposition it suffices to show that for some $\beta>0$
which depends only on $p$, we have
%
%
\begin{equation}
\label{l21} P\bigl(Y_{i+1}-Y_i\geq m\vert F\cap
G(x_0,\ldots, x_i,r_1,\ldots,r_i)
\bigr)\leq\exp(-\beta m).
\end{equation}
By Lemma \ref{l10} on the event $F\cap G(x_0,\ldots, x_i,r_1,\ldots
,r_i)\cap\{Y_{i+1}-Y_i\geq m\}$, there is a point in $u\in D_i$ such that
$u\rightarrow(\bigcup_{j=r_i+m}^n L_j)\cap C_{x_i,r_i}$.
Hence, \eqref{l21} will follow if we prove
%
%
\begin{eqnarray}
\label{l22} %
&&P\Biggl(u\rightarrow C_{x_i,r_i}\cap\Biggl(\bigcup
_{j=r_i+m}^n L_j\Biggr)\mbox{ for
some } u \in D_i\Big\vert
\nonumber
\\[-8pt]
\\[-8pt]
\nonumber
&& \hspace*{85pt}F\cap G(x_0,\ldots, x_i,r_1,
\ldots,r_i)\Biggr)\leq\exp(-\beta m).
\end{eqnarray}
By Lemma \ref{l10} this can be written as
\begin{eqnarray*}
&&P \Biggl(u\rightarrow C_{x_i,r_i}\cap\Biggl(\bigcup
_{j=r_i+m}^n L_j\Biggr)\mbox{ for some } u
\in D_i\Big\vert
\\
&&\hspace*{18pt}F\cap\bigl\{(x_i,r_i)
\rightarrow L_n\bigr\}\cap\{u\nrightarrow L_n\ \forall u\in
D_i\}
\cap\{u \nrightarrow L_{r_i+1}\ \forall u\in
B_i\}\\
&&\hspace*{108pt}{}\cap\Biggl(\bigcap_{j=1}^i
\{u\nrightarrow L_n,u\nrightarrow C_{x_{j-1},r_{j-1}}\ \forall u\in
A_{i,j}\} \Biggr) \Biggr)
\\
&&\qquad\leq\exp(-\beta m).
\end{eqnarray*}
Since the state of the bonds above $L_{r_i}$ is independent of $F$ this
is equivalent to
\begin{eqnarray*}
&&P \Biggl(u\rightarrow C_{x_i,r_i}\cap\Biggl(\bigcup
_{j=r_i+m}^n L_j\Biggr)\mbox{ for some } u
\in D_i \Big\vert
\\
&&\hspace*{15pt}\bigl\{(x_i,r_i)
\rightarrow L_n\bigr\}\cap\{u\nrightarrow L_n\ \forall u\in
D_i\}
\cap\{u \nrightarrow L_{r_i+1}\ \forall u\in
B_i\}\\
&&\hspace*{84pt}{}\cap\Biggl(\bigcap_{j=1}^i
\{u\nrightarrow L_n,u\nrightarrow C_{x_{j-1},r_{j-1}}\ \forall u\in
A_{i,j}\} \Biggr) \Biggr)
\\
&&\qquad\leq\exp(-\beta m).
\end{eqnarray*}
%
Since $\{(x_i,r_i)\rightarrow L_n\}$ is a disjoint union of the events
$\{\Gamma_r=\gamma\}$ where $\gamma$ ranges over $\Upsilon$, it
suffices to show that for
any $\gamma\in\Upsilon$ we have
%
%
\begin{eqnarray}
\label{eq23} %
&& P \Biggl(u\rightarrow C_{x_i,r_i}\cap\Biggl(
\bigcup_{j=r_i+m}^n L_j\Biggr)
\mbox{ for some } u \in D_i \Big\vert
\nonumber\\
&&\hspace*{18pt} \{\Gamma_r=\gamma
\}\cap\{u\nrightarrow L_n\ \forall u\in D_i\}
\cap\{u
\nrightarrow L_{r_i+1}, \forall u\in B_i\}
\nonumber
\\[-8pt]
\\[-8pt]
\nonumber
&&\hspace*{61pt}{}\cap\Biggl(\bigcap
_{j=1}^i\{u\nrightarrow L_n,u
\nrightarrow C_{x_{j-1},r_{j-1}}\ \forall u\in A_{i,j}\} \Biggr)
\Biggr)
\\
&&\qquad\leq\exp(-\beta m).\nonumber
\end{eqnarray}
But, as explained at the beginning of this proof, the left-hand side
above can be written as
%
%
\begin{eqnarray}
\label{eq24} %
&& P \Biggl(u\rightarrow C_{x_i,r_i}\cap\Biggl(
\bigcup_{j=r_i+m}^n L_j\Biggr)
\mbox{ for some } u \in D_i \Big\vert\nonumber\\
&&\hspace*{15pt} \{\gamma\mbox{ is open}\}\cap
D(\gamma)\cap\{u\nrightarrow L_n\ \forall u\in D_i\}\cap
\{u \nrightarrow L_{r_i+1}\ \forall u\in B_i\}
\\
&&\hspace*{102pt}{}\cap\Biggl(
\bigcap_{j=1}^i\{u\nrightarrow
L_n,u\nrightarrow C_{x_{j-1},r_{j-1}}\ \forall u\in A_{i,j}\}
\Biggr) \Biggr).\nonumber
\end{eqnarray}
Now, let $V(\gamma)$ be the set of vertices of $\gamma$. Then, noting that
%
%
\begin{eqnarray}
&& \Biggl\{u\rightarrow C_{x_i,r_i}\cap\Biggl(\bigcup
_{j=r_i+m}^n L_j\Biggr)\mbox{ for some } u
\in D_i\Biggr\}\nonumber\\
&&\quad{}\cap
 \{u\nrightarrow L_n\ \forall u\in
D_i\}\cap\{\gamma\mbox{ is open}\}
\nonumber
\\[-8pt]
\\[-8pt]
\nonumber
&&\qquad =\Biggl\{u\rightarrow
C_{x_i,r_i}\cap\Biggl(\bigcup_{j=r_i+m}^n
L_j\Biggr)\mbox{ within } A_{\gamma, r}\mbox{ for some } u \in
D_i\Biggr\}
\\
&&\qquad\quad{} \cap\{u\nrightarrow L_n\ \forall u\in
D_i\} \cap\{\gamma\mbox{ is open}\},\nonumber
\end{eqnarray}
and that
%
%
\begin{eqnarray}
&& \{u\nrightarrow L_n\ \forall u\in D_i\}\cap
\{\gamma\mbox{ is open}\}
\nonumber
\\[-8pt]
\\[-8pt]
\nonumber
&&\qquad=\bigl\{u\nrightarrow L_n\cup V(\gamma)\
\forall u\in D_i\bigr\}\cap\{\gamma\mbox{ is open}\}, %
\end{eqnarray}
\eqref{eq24} can be written as
%
%
\begin{eqnarray}
\label{eq25} %
&& P\Biggl(u\rightarrow C_{x_i,r_i}\cap\Biggl(
\bigcup_{j=r_i+m}^n L_j\Biggr)
\mbox{ within } A_{\gamma, r}\mbox{ for some } u \in D_i\Big\vert
\nonumber\\
&&\hspace*{15pt}
\{\gamma\mbox{ is open}\}\cap D(\gamma)\cap\bigl\{u\nrightarrow L_n
\cup V(\gamma)
\ \forall u\in D_i\bigr\}
\\
&&{}\hspace*{12pt} \cap\{u \nrightarrow
L_{r_i+1}\ \forall u\in B_i\}\nonumber \\
&&\hspace*{64pt}{}\cap\Biggl(\bigcap
_{j=1}^i\{u\nrightarrow L_n,u
\nrightarrow C_{x_{j-1},r_{j-1}}\ \forall u\in A_{i,j}\}
\Biggr)\Biggr).\nonumber
\end{eqnarray}
Since $\{\gamma\mbox{ is open}\}$ is independent of all the other
events involved in the above expression, \eqref{eq25} is equal to
%
%
\begin{eqnarray}
\label{eq26} %
& &P\Biggl(u\rightarrow C_{x_i,r_i}\cap\Biggl(
\bigcup_{j=r_i+m}^n L_j\Biggr)
\mbox{ within } A_{\gamma, r}\mbox{ for some } u \in D_i \Big\vert
\nonumber\\
&&\hspace*{16pt} D(\gamma)\cap\bigl\{u\nrightarrow L_n \cup V(\gamma)\ \forall
u\in
D_i\bigr\} \cap\{u \nrightarrow L_{r_i+1}\ \forall u\in
B_i\}\\
&&\hspace*{79pt}{} \cap\Biggl(\bigcap_{j=1}^i
\{u\nrightarrow L_n,u\nrightarrow C_{x_{j-1},r_{j-1}}\ \forall u\in
A_{i,j}\}\Biggr)\Biggr).\nonumber
\end{eqnarray}
Since $\{u\rightarrow C_{x_i,r_i}\cap(\bigcup_{j=r_i+m}^n L_j)\mbox{
within } A_{\gamma, r}\mbox{ for some } u \in D_i\}$
is an increasing event while
$ D(\gamma)\cap\{u\nrightarrow L_n\cup V(\gamma)\ \forall u\in D_i\}
\cap\{u \nrightarrow L_{r_i+1}\ \forall u\in B_i\}\cap(\bigcap
_{j=1}^i\{u\nrightarrow L_n,u\nrightarrow C_{x_{j-1},r_{j-1}} \ \forall
u\in A_i\} ) $
is decreasing, \eqref{eq26}
is bounded above by
%
%
\begin{eqnarray}
&&P\Biggl(u\rightarrow C_{x_i,r_i}\cap\Biggl(\bigcup
_{j=r_i+m}^n L_j\Biggr)\mbox{ within }
A_{\gamma, r}\mbox{ for some } u \in D_i\Biggr)
\nonumber
\\
&&\qquad \leq P\Biggl(u
\rightarrow C_{x_i,r_i}\cap\Biggl(\bigcup_{j=r_i+m}^n
L_j\Biggr) \mbox{ for some } u\in L_{r_i}\\
&&\hspace*{148pt}\mbox{to the
right of } (x_i,r_i)\Biggr).\nonumber
\end{eqnarray}
The proposition now follows from Lemma \ref{l0}.
\end{pf}

\begin{pf*}{Proof of Lemma \ref{c1}}
If $I\geq m$ then the sequence $Y_0,\ldots,Y_{m}$ is strictly
increasing. But it follows from Proposition~\ref{p1} that
$P(Y_{i+1}>Y_i\vert Y_1,\ldots,Y_i)\leq\exp(-\beta)$ a.s. Hence
(a) follows by induction in $i$. To prove (b)
write
\begin{eqnarray*}
P\bigl(Y_I\geq m^2\bigr)&\leq& P(I>m)+P
\bigl(Y_m\geq m^2\bigr)
\\
&\leq&\exp(-\beta m) +\sum_{j=0}^{m-1}
P(Y_{j+1}-Y_j\geq m)
\\
&\leq&\exp(-\beta m)+m \exp(-\beta m),
\end{eqnarray*}
where the second inequality follows from part (a) and the last one from
Proposition~\ref{p1}. Part (c) follows easily from parts (a) and (b).
\end{pf*}
%




%
\begin{lemma}\label{c5}
Let $A'\in S\cup S_0$ and let $f$ be a cylinder function on $S\cup S_0$
depending only on the coordinates $-2r,\ldots,-2$. Then,
$|E(f(\zeta^{A'}_n)\vert\{0\rightarrow L_n\})-E(f(\zeta^{0}_n)\vert\{
0\rightarrow L_n\})|\leq2(n+r)\|f\|\exp(-\beta n)\ \forall n \in\N$.
\end{lemma}
\begin{pf}
Let $\Phi$ be the (random) set of points in levels $1,2,\ldots, n$ which
can be reached from $(0,0)$ following an open path.
The event $\{0\rightarrow L_n\}$ is a disjoint union of events of the
form $\{\Phi= \kappa\}$ where $\kappa$ ranges over all values of
$\Phi$ containing at least one point in
$L_n$. Then write
\begin{eqnarray*}
&&\bigl |E\bigl(f\bigl(\zeta^{A'}_n\bigr)\vert\Phi= \kappa\bigr)-E
\bigl(f\bigl(\zeta^{0}_n\bigr)\vert\Phi=\kappa\bigr)\bigr|
\\
&&\qquad\leq2\|f\| \sum_{i=1}^{n+r}P(-2i\rightarrow
L_n \mbox{ off }\kappa)
\\
&&\qquad\leq2(n+r)\|f\|\exp(-\beta n)
\end{eqnarray*}
and the lemma follows.
\end{pf}

Before starting the proof of Lemma \ref{c6}, we need to introduce some
further notation: $T$ will the map sending subsets of $2\Z+k$ into
subests of $2\Z+k+1$ given by
$T(A)=\{x-1;x\in A\}$ and for $(x,n)\in\Lambda$ let $\mathcal
{G}_{x,n}^+$ be the $\sigma$-algebra generated by the random variables
determining the state of the bonds whose vertices
are in $\bigcup_{i=0}^n L_i$ and by the bonds having at least one vertex
strictly to the right of $\{(x+i,n+i); i=0,1,\ldots\}$ and let
$\mathcal
{G}_{x,n}^-$ be the $\sigma$-algebra generated by the random variables
determining the state of all the
other bonds.
If $B$ is an infinite subset of $2\Z+k$ which is bounded above, we
define for $n\geq k$:
$\xi_{k,n}^B=\{z\dvtx(x,k)\rightarrow(z,n) \mbox{ for some } x\in
B\}$,
$r(\xi_{k,n}^B)= \sup(\xi_{k,n}^B)$ and
$\zeta^B_{k,n}= \{x-r(\xi^B_{k,n})\dvtx x \in\xi^B_{k,n}\}$. As before
$(\zeta^B_{k,n},n\geq k) $ is a Markov chain
on infinite subsets of $2\Z_-$ containing $0$.


\begin{pf*}{Proof of Lemma \ref{c6}}
Since $E(f(\zeta_n^A)\vert I=i, Y_i=j)$ is a convex combination of
$E(f(\zeta_n^A)\vert I=i, Y_i=j, X_i=x_i)$ where $ x_i$ runs over all
possible values of~$X_i$, it suffices to show that for all $x_i$
%
%
\begin{eqnarray}
\label{f1} %
&&\bigl\vert E\bigl(f\bigl(\zeta_n^A
\bigr)\vert I=i, Y_i=j, X_i=x_i\bigr)-E\bigl(f
\bigl(\zeta^0_n\bigr)\vert0\rightarrow L_n
\bigr)\bigr\vert
\nonumber
\\[-8pt]
\\[-8pt]
\nonumber
&&\qquad\leq2\|f\| (n+r)\exp\bigl(-\beta(n-j)\bigr).
\end{eqnarray}
But on the event $\{I=i,Y_i=j,X_i=x_i\}$ it happens that $x_i$ is the
rightmost point of $\xi_j^A$ from which there is an open path to $L_n$.
Therefore,
\eqref{f1} will follow if we show that for all infinite subset $A'$ of
$L_j$ such that $\sup A'=x_i$ we have
%
%
\begin{eqnarray}
\label{f2} %
&&\bigl\vert E\bigl(f\bigl(\zeta_{j,n}^{A'}
\bigr)\vert I=i, Y_i=j, X_i=x_i\bigr)-E\bigl(f
\bigl(\zeta^0_n\bigr)\vert0\rightarrow L_n
\bigr)\bigr\vert
\nonumber
\\[-8pt]
\\[-8pt]
\nonumber
&&\qquad\leq2\|f\| (n+r)\exp\bigl(-\beta(n-j)\bigr).
\end{eqnarray}
%
Since $\{I=i,Y_i=j,X_i=x_i\}=\{(x_i,j)\rightarrow L_n\}\cap H$ where
$H\in{\mathcal G}^+_{x_i,j}$ and the evolution of $\zeta^{A'}_{j,k}$
as $k$ increases
from $j$ to $n$ is
${\mathcal G}^-_{x_i,j}$-measurable we have
%
%
\begin{equation}
\label{f3} E\bigl(f\bigl(\zeta_{j,n}^{A'}\bigr)\vert I=i,
Y_i=j, X_i=x_i\bigr)=E\bigl(f\bigl(\zeta
_{j,n}^{A'}\bigr)\vert(x_i,j)\rightarrow
L_n\bigr).
\end{equation}
Hence, \eqref{f2} follows from Lemma \ref{c5} and translation invariance.
\end{pf*}

\section*{Acknowledgements}
The author wishes to thank A. Galves who mentioned to him the
plausibility of Theorem \ref{T1}, P. Ferrari for conversations
concerning the results of
\cite{FKM}, F. Ezanno who provided Figure \ref{fig1} and careful reading which
improved this paper, and the referees for their useful comments.

%

%



\printaddresses

\end{document}